\documentstyle[12pt]{article}
         \title{On $L^{1}$-Convergence of
Fourier Series Under $MVBV$ Condition}
\author{Dansheng Yu,\thanks{Supported in part
by NSERC RCD grant of St. Francis Xavier University and in part by AARMS of Canada.}
 Ping Zhou\thanks{Supported by NSERC of Canada.} and Songping Zhou\thanks{The third author's research is done as a W. F. James Chair Professor
of St. Francis Xavier University. His research is also supported in
part by NSF of China under grant number 10471130.}}

\date{}

\begin{document}
          \maketitle
          \pagenumbering{arabic}
\begin{abstract}
{\footnotesize Let $f\in L_{2\pi }$ be a real-valued even function with its Fourier series $%
\frac{a_{0}}{2}+\sum_{n=1}^{\infty }a_{n}\cos nx,$ and let
$S_{n}\left( f,x\right) ,\;n\geq 1,$ be the $n$-th partial sum of
the Fourier series. It is well-known that if the nonnegative
sequence $\{a_{n}\}$ is decreasing and
$\lim\limits_{n\rightarrow \infty }a_{n}=0$, then%
$$
\lim\limits_{n\rightarrow \infty }\Vert f-S_{n}(f)\Vert _{L}=0\;\;%
\mbox{if and only if}\;\;\lim\limits_{n\rightarrow \infty }a_{n}\log
n=0.
$$
We weaken the monotone condition in this classical result to the
so-called mean value bounded variation ($MVBV$) condition. The
generalization of the above classical result in real-valued function
space is presented as a
special case of the main result in this paper which gives the $L^{1}$%
-convergence of a function $f\in L_{2\pi }$ in complex space. We
also give
results on $L^{1}$-approximation of a function $f\in L_{2\pi }$ under the $%
MVBV$ condition.}
\end{abstract}

\begin{quote}



\end{quote}

\normalsize
\section{Introduction}

Let $L_{2\pi }$ be the space of all complex-valued integrable functions $%
f(x) $ of period $2\pi $ equipped with the norm
$$
\Vert f\Vert _{L}=\int_{-\pi }^{\pi }|f(x)|dx.$$
 Denote the Fourier
series of $f\in L_{2\pi }$ by
\begin{eqnarray}
\sum\limits_{k=-\infty }^{\infty }\hat{f}(k)e^{ikx},  \label{1}
\end{eqnarray}%
and its partial sum $S_{n}(f,x)$ by
$$
\sum\limits_{k=-n}^{n}\hat{f}(k)e^{ikx}.
$$
When $f(x)\in L_{2\pi }$ is a real valued even function, then the
Fourier series of $f$ has the form
\begin{eqnarray}
\frac{a_{0}}{2}+\sum\limits_{k=1}^{\infty }a_{k}\cos kx,  \label{1a}
\end{eqnarray}%
correspondingly, its partial sum $S_{n}(f,x)$ is
$$
\frac{a_{0}}{2}+\sum\limits_{k=1}^{n}a_{k}\cos kx.
$$

The following two classical convergence results can be found in many
monographs (see \cite{Boas} and \cite{Zygmund}, for example):

\textbf{Result One:} If a nonnegative sequence
$\{b_{n}\}_{n=1}^{\infty }$
is decreasing and $\lim_{n\rightarrow \infty }b_{n}=0$, then the series $%
\sum\limits_{n=1}^{\infty }b_{n}\sin nx$ converges uniformly if and only if $%
\lim\limits_{n\rightarrow \infty }nb_{n}=0.$

\textbf{Result Two:} Let $f\in L_{2\pi }$ be an even function and
(\ref{1a}) be its Fourier series. If the sequence
$\{a_{n}\}_{n=0}^{\infty }$ is nonnegative, decreasing, and
$\lim_{n\rightarrow \infty }a_{n}=0$, then
$$
\lim\limits_{n\rightarrow \infty }\Vert f-S_{n}(f)\Vert _{L}=0\;\;%
\mbox{if and only if}\;\;\lim\limits_{n\rightarrow \infty }a_{n}\log
n=0.
$$

These results have been generalized by weakening the monotone
conditions of the coefficient sequences. They have also been
generalized to the complex valued function spaces. The most recent
generalizations of Result One can be found in \cite{Zhou-Zhou-Yu}
where the monotonic condition is finally weakened to the $MVBV$
condition (Mean Value Bounded Variation condition, see Corollary 2
in Section 2 for definition), and it is proved to be the weakest
possible condition we can have to replace the monotone condition in
Result One. The process of generalizing Result Two can be found in
many papers, for example, see \cite{Le-Zhou} - \cite{Yu-Zhou}. In
this paper, we will weaken the monotone condition in Result Two (and
all its later generalized conditions, see \cite{Zhou-Zhou-Yu} for
the relations between these conditions), to the $MVBV$ condition in
the complex valued function spaces (see \textit{Definition 1 }in
Section 2) in Theorem 1, and give the generalization in real valued
function spaces as a special case of Theorem 1 in Corollary 2. Like
the important role that the $MVBV$ condition plays in generalizing
Result One, although we are not able to prove it here, we propose
that Theorem 1 in Section 2 is the ultimate generalization of Result
Two, i.e. the $MVBV$ condition is also the weakest possible
condition we can have to replace the monotone condition in Result
Two. We also discuss, under the $MVBV$ condition, the
$L^{1}$-approximation rate of a function $f\in L_{2\pi }$ in the
last section.

Throughout this paper, we always use $C(x)$ to indicate a positive
constant depending upon $x$ only, and use $C$ to indicates an
absolute positive constant. They may have different values in
different occurrences.

\section{$L^{1}$ convergence}

In this section, we first give the definition of $%
MVBV$ condition, or the class $MVBVS$, and then prove our main result on $%
L^{1}$-convergence of the Fourier series of a\ complex valued function $%
f(x)\in L_{2\pi }$ whose coefficients form a sequence in the class
$MVBVS$.

{\bf Definition 1.}\quad {\it Let
$\mathbf{c}:=\{c_{n}\}_{n=0}^{\infty }$ be a sequence of complex
numbers satisfying $c_{n}\in K(\theta _{1}):=\{z:|\arg z|\leq \theta
_{1}\}$ for some $\theta _{1}\in \lbrack 0,\pi /2)$ and all
$n=0,1,2,\cdots .$ If there is a number $\lambda \geq 2$ such that
$$
\sum\limits_{k=m}^{2m}|\Delta
c_{k}|:=\sum\limits_{k=m}^{2m}|c_{k+1}-c_{k}|\leq C(\mathbf{c})\frac{1}{m}%
\sum\limits_{k=[\lambda ^{-1}m]}^{[\lambda m]}|c_{k}|
$$
holds for all $m=1,2,\cdots ,$ then we say that the sequence
$\mathbf{c}$ is a Mean Value Bounded Variation Sequence, i.e.,
$\mathbf{c}\in MVBVS,$ in complex sense, or the sequence
$\mathbf{c}$ satisfies the $MVBV$ condition.}

Our main result of this paper is:

{\bf Theorem 1.}\quad {\it Let $f(x)\in L_{2\pi }$ be a
complex-valued function. If the Fourier coefficients $\hat{f}(n)$ of
$f$ satisfy that $\{\hat{f}(n)\}_{n=0}^{+\infty }\in MVBVS$ and
\begin{eqnarray}
\lim\limits_{\mu \rightarrow 1^{+}}\limsup\limits_{n\rightarrow
\infty }\sum\limits_{k=n}^{[\mu n]}|\Delta \hat{f}(k)-\Delta
\hat{f}(-k)|\log k=0, \label{2}
\end{eqnarray}%
where%
$$
\Delta \hat{f}(k)=\hat{f}(k+1)-\hat{f}(k),\;\;\Delta \hat{f}(-k)=\hat{f}%
(-k-1)-\hat{f}(-k),\;\;k\geq 0.
$$
Then
$$
\lim\limits_{n\rightarrow \infty }\Vert f-S_{n}(f)\Vert _{L}=0
$$
if and only if
$$
\lim\limits_{n\rightarrow \infty }\hat{f}(n)\log |n|=0.
$$
}

 In order to prove Theorem 1, we present the
following four lemmas.

{\bf Lemma 1.}\quad {\it Let $\{c_{n}\}\in MVBVS,$ then for any
given $1<\mu <2$, we have
$$
\sum\limits_{k=n}^{[\mu n]}|\Delta c_{k}|\log k=O\left(
\max\limits_{[\lambda ^{-1}n]\leq k\leq \lbrack \lambda
n]}|c_{k}|\log k\right) ,\;\;\;\;n\rightarrow \infty ,
$$
where the implicit constant depends only on the sequence $\{c_{n}\}$ and $%
\lambda .$ }

For sufficiently large $n$, the lemma can be derived directly from
the conditions that $1<\mu <2$ and $\{c_{n}\}\in MVBVS.$

{\bf Lemma 2.}\quad {\it Let $\{\hat{f}(n)\}\in K(\theta _{0})$ for
some $\theta _{0}\in \lbrack 0,\pi /2)$, then
$$
\sum\limits_{k=1}^{n}\frac{1}{k}\left| \hat{f}(n+k)\right| =O\left(
\Vert f-S_{n}(f)\Vert _{L}\right)
$$
for all $n=1,2,\cdots ,$ where the implicit constant depends only on
$\theta _{0}.$ }

{\bf Proof.} Write
$$
\phi _{\pm n}(x):=\sum\limits_{k=1}^{n}\frac{1}{k}\left( e^{i(k\mp
n)x}-e^{-i(k\pm n)x}\right) .
$$
It follows from a well-known inequality (e.g. see Theorem 2.5 in \cite%
{Xie-Zhou2})
$$
\sup\limits_{n\geq 1}\left| \sum\limits_{k=1}^{n}\frac{\sin
kx}{k}\right| \leq 3\sqrt{\pi }
$$
that
$$
|\phi _{\pm n}(x)|\leq 6\sqrt{\pi }.
$$
Hence
$$
\frac{1}{6\sqrt{\pi }}\left| \int_{-\pi }^{\pi
}(f(x)-S_{n}(f,x))\phi _{\pm n}(x)dx\right| \leq \Vert
f-S_{n}(f)\Vert _{L},
$$
and therefore
$$
\left| \sum\limits_{k=1}^{n}\frac{1}{k}\hat{f}(n+k)\right| =O(\Vert
f-S_{n}(f)\Vert _{L}).
$$
Now as $\{\hat{f}(n)\}\in K(\theta _{0})$ for some $\theta _{0}\in
\lbrack
0,\pi /2)$ and for all $n\geq 1,$ we have%
\begin{eqnarray*}
\sum\limits_{k=1}^{n}\frac{1}{k}\left| \hat{f}(n+k)\right| &\leq
&C(\theta_{0})\sum\limits_{k=1}^{n}\frac{1}{k}\mbox{Re}\hat{f}(n+k) \\
&\leq &C(\theta _{0})\left| \sum\limits_{k=1}^{n}\frac{1}{k}\hat{f}%
(n+k)\right| \\
&=&O(\Vert f-S_{n}(f)\Vert _{L}).
\end{eqnarray*}

{\bf Lemma 3.}\quad  (\cite{Xie-Zhou}). {\it Write
\begin{eqnarray*}
D_{k}(x) &:&=\frac{\sin ((2k+1)x/2)}{2\sin (x/2)}, \\
&&\;\;\;\;\ \ \  \\
D_{k}^{\ast }(x) &:&=\left\{
\begin{array}{ll}
\frac{\cos (x/2)-\cos ((2k+1)x/2)}{2\sin (x/2)} & |x|\leq 1/n, \\
-\frac{\cos ((2k+1)x/2)}{2\sin (x/2)} & 1/n\leq |x|\leq \pi ,%
\end{array}%
\right. \\
&&\;\;\;\;\ \ \  \\
E_{k}(x) &:&=D_{k}(x)+iD_{k}^{\ast }(x).
\end{eqnarray*}
For $k=n,n+1,\cdots ,2n$, we have
\begin{eqnarray}
E_{k}(\pm x)-E_{k-1}(\pm x)=e^{\pm ikx},  \label{3}
\end{eqnarray}%
\begin{eqnarray}
E_{k}(x)+E_{k}(-x)=2D_{k}(x),  \label{4}
\end{eqnarray}%
\begin{eqnarray}
\Vert E_{k}\Vert _{L}+\Vert D_{k}\Vert _{L}=O(\log k).  \label{5}
\end{eqnarray}
}

{\bf Lemma 4.}\quad {\it Let $\{\hat{f}(n)\}\in MVBVS$. If
$\lim\limits_{n\rightarrow \infty }\Vert f-S_{n}(f)\Vert _{L}=0,$
then
$$
\lim\limits_{n\rightarrow \infty }\hat{f}(n)\log n=0.
$$}

{\bf Proof.} By the definition of $MVBVS$, we derive that for
$k=n,n+1,\cdots ,2n,$ \begin{eqnarray*} |\hat{f}(2n)| &\leq
&\sum\limits_{j=k}^{2n-1}|\Delta \hat{f}(j)|+|\hat{f}(k)|
\\
&\leq &\sum\limits_{j=k}^{2k}|\Delta \hat{f}(j)|+|\hat{f}(k)| \\
&=&O\left( \frac{1}{n}\sum\limits_{j=[\lambda ^{-1}k]}^{[\lambda k]}|\hat{f}%
(j)|\right) +|\hat{f}(k)|.
\end{eqnarray*}
Therefore, it follows that from the fact that%
$$
\log n\leq C\left( \lambda \right) \sum\limits_{j=[\lambda
]+1}^{[(\lambda +1)^{-2}n]}\frac{1}{j},
$$
we have%
\begin{eqnarray}
|\hat{f}(2n)|\log n &\leq &C(\lambda
)|\hat{f}(2n)|\sum\limits_{j=[\lambda
]+1}^{[(\lambda +1)^{-2}n]}\frac{1}{j}  \nonumber \\
&\leq &C(\lambda )\sum\limits_{j=[\lambda ]+1}^{[(\lambda +1)^{-2}n]}\frac{1%
}{j}\left( \frac{1}{n}\sum\limits_{k=[\lambda ^{-1}(n+j)]}^{[\lambda (n+j)]}|%
\hat{f}(k)|+|\hat{f}(n+j)|\right)  \nonumber \\
&=&\frac{C(\lambda )}{n}\sum\limits_{j=[\lambda ]+1}^{[(\lambda +1)^{-2}n]}%
\frac{1}{j}\sum\limits_{k=[\lambda ^{-1}(n+j)]}^{[\lambda
(n+j)]}|\hat{f}(k)|
\nonumber\\
&&\;\;\;\;\;\;+C(\lambda )\sum\limits_{j=1}^{[(\lambda +1)^{-2}n]}\frac{1}{j}%
|\hat{f}(n+j)|  \nonumber \\
&=&:I_{1}+I_{2}.  \label{6}
\end{eqnarray}%
By applying Lemma 2, we see that
\begin{eqnarray}
I_{2}\leq C(\lambda ,\theta _{0})\Vert f-S_{n}(f)\Vert _{L}.
\label{7}
\end{eqnarray}%
We calculate $I_{1}$ as follows (note that we may add more repeated
terms in the right hand side of every inequality below):
\begin{eqnarray}
I_{1} &\leq &\frac{C(\lambda )}{n}\sum\limits_{j=[\lambda
]+1}^{[(\lambda +1)^{-2}n]}\frac{1}{j}\sum\limits_{k=[\lambda
^{-1}n]+[\lambda
^{-1}j]}^{[\lambda n]+[\lambda j]+1}|\hat{f}(k)|  \nonumber\\
&\leq &\frac{C(\lambda )}{n}\sum\limits_{j=[\lambda ]+1}^{[(\lambda
+1)^{-2}n]}\sum\limits_{m=1}^{[(\lambda
+1)^{2}]}\sum\limits_{k=[\lambda ^{-1}n]}^{[\lambda
n]+1}\frac{\left| \hat{f}\left( m[\lambda
^{-1}j]+k\right) \right| }{j}  \nonumber\\
&\leq &\frac{C(\lambda )}{n}\sum\limits_{m=1}^{[(\lambda
+1)^{2}]}\sum\limits_{j=[\lambda ]+1}^{[(\lambda
+1)^{-2}n]}\sum\limits_{k=0}^{[\lambda n]-[\lambda
^{-1}n]+1}\frac{\left| \hat{f}\left( [\lambda ^{-1}n]+m[\lambda
^{-1}j]+k\right) \right| }{j}
\nonumber\\
&\leq &\frac{C(\lambda )}{n}\sum\limits_{m=1}^{[(\lambda
+1)^{2}]}\sum\limits_{k=0}^{[\lambda n]-[\lambda
^{-1}n]+1}\sum\limits_{j=1}^{m\left[ \left( \lambda (\lambda
+1)^{2}\right) ^{-1}n\right] }\frac{\left| \hat{f}\left( [\lambda
^{-1}n]+k+j\right)
\right| }{j}  \nonumber\\
&&\;\;\;\;  \nonumber
\\
&\leq &\frac{C(\lambda
)}{n}\sum\limits_{m=1}^{[(\lambda+1)^{2}]}\sum\limits_{k=0}^{[\lambda
n]-[\lambda ^{-1}n]+1}\left\|
f-S_{[\lambda ^{-1}n]+k}(f)\right\| _{L}\;\;\mbox{(by Lemma 2)}  \nonumber\\
&\leq &\frac{C(\lambda )}{n}\sum\limits_{k=0}^{[\lambda n]-[\lambda
^{-1}n]+1}\left\| f-S_{[\lambda ^{-1}n]+k}(f)\right\| _{L}.
\label{8}
\end{eqnarray}%
Finally, by combining (\ref{6}) - (\ref{8}) and the condition
$$
\lim\limits_{n\rightarrow \infty }\Vert f-S_{n}(f)\Vert _{L}=0,
$$
we get
$$
\lim\limits_{n\rightarrow \infty }\hat{f}(2n)\log n=0.
$$
A similar argument yields that
$$
\lim\limits_{n\rightarrow \infty }|\hat{f}(2n+1)|\log n=0.
$$
This proves Lemma 4.

 We now come to the proof of Theorem 1.

{\bf Proof of Theorem 1.}
\textbf{Sufficiency.} Given $\varepsilon >0$, by (\ref{2}), there is a $%
1<\mu <2$ such that
\begin{eqnarray}
\sum\limits_{k=n}^{[\mu n]}\left| \Delta \hat{f}(k)-\Delta \hat{f}%
(-k)\right| \log k\leq \varepsilon  \label{9}
\end{eqnarray}%
holds for sufficiently large $n>0$. Let
$$
\tau _{\mu n,n}(f,x):=\frac{1}{[\mu n]-n}\sum\limits_{k=n}^{[\mu
n]-1}S_{k}(f,x)
$$
be the Vall\'{e}e Poussin sum of order $n$ of $f$. Then we have
\begin{eqnarray}
\lim\limits_{n\rightarrow \infty }\Vert f-\tau _{\mu n,n}(f)\Vert
_{L}=0. \label{10}
\end{eqnarray}%
By (\ref{3}), (\ref{4}), and applying Abel transformation, we get
\begin{eqnarray}
&&\tau _{\mu n,n}(f,x)-S_{n}(f,x)\nonumber \\
&=&\frac{1}{[\mu n]-n}\sum\limits_{k=n+1}^{[\mu n]}([\mu n]-k)\left( \hat{f}%
(k)e^{ikx}+\hat{f}(-k)e^{-ikx}\right)  \nonumber\\
&=&\frac{1}{[\mu n]-n}\sum\limits_{k=n}^{[\mu n]}([\mu n]-k)\left(
2\Delta \hat{f}(k)D_{k}(x)-(\Delta \hat{f}(k)-\Delta
\hat{f}(-k))E_{k}(-x)\right)
\nonumber\\
&&+\frac{1}{[\mu n]-n}\sum\limits_{k=n}^{[\mu n]-1}\left( \hat{f}%
(k+1)E_{k}(x)-\hat{f}(-k-1)E_{k}(-x)\right)  \nonumber\\
&&-\left( \hat{f}(n)E_{n}(x)+\hat{f}(-n)E_{n}(-x)\right) .
\label{11}
\end{eqnarray}%
Thus, by (\ref{5}) and Lemma 1, we have
\begin{eqnarray}
&&\Vert f-S_{n}(f)\Vert _{L}  \nonumber\\
&\leq &\Vert f-\tau _{\mu n,n}(f)\Vert _{L}+\Vert \tau _{\mu
n,n}(f)-S_{n}(f)\Vert _{L}  \nonumber\\
&=&\Vert f-\tau _{\mu n,n}(f)\Vert _{L}+O\left(
\sum\limits_{k=n}^{[\mu
n]}\left| \Delta \hat{f}(k)\right| \log k\right)  \nonumber \\
&&+O\left( \sum\limits_{k=n}^{[\mu n]}\left| \Delta \hat{f}(k)-\Delta \hat{f}%
(-k)\right| \log k\right)  \nonumber\nonumber \\
&&+O\left( \max\limits_{n\leq |k|\leq \lbrack \mu
n]}|\hat{f}(k)|\log
|k|\right)  \nonumber\\
&=&\Vert f-\tau _{\mu n,n}(f)\Vert _{L}+O\left(
\max\limits_{[\lambda
^{-1}n]\leq |k|\leq \lbrack \lambda n]}|\hat{f}(k)|\log |k|\right)  \nonumber\\
&&+O\left( \sum\limits_{k=n}^{[\mu n]}\left| \Delta \hat{f}(k)-\Delta \hat{f}%
(-k)\right| \log k\right) ,  \label{11a}
\end{eqnarray}%
then
$$
\limsup\limits_{n\rightarrow \infty }\Vert f-S_{n}(f)\Vert _{L}\leq
\varepsilon
$$
follows from (\ref{9}), (\ref{10}) and the condition that
$$
\lim\limits_{n\rightarrow \infty }\hat{f}(n)\log |n|=0.
$$
This implies that
$$
\lim\limits_{n\rightarrow \infty }\Vert f-S_{n}(f)\Vert _{L}=0.
$$

\textbf{Necessity.} Since $\{\hat{f}(n)\}\in MVBVS$, by applying
Lemma 4, we have
\begin{eqnarray}
\lim\limits_{n\rightarrow \infty }\hat{f}(n)\log n=0.  \label{12}
\end{eqnarray}%
In order to prove $\lim\limits_{n\rightarrow -\infty }\hat{f}(n)\log
\left| n\right| =0,$ by applying (\ref{11}) and (\ref{5}), we see
that for any given $\mu ,$ $1<\mu <2,$
\begin{eqnarray}
&&\Vert \hat{f}(-n)E_{n}(-x)\Vert _{L} \nonumber\\
&\leq &\Vert \tau _{\mu n,n}(f)-S_{n}(f)\Vert _{L}  \nonumber \\
&&+\frac{1}{[\mu n]-n}\left\| \sum\limits_{k=n}^{[\mu n]-1}\hat{f}%
(-k-1)E_{k}(-x)\right\| _{L}  \nonumber\\
&&+O\left( \sum\limits_{k=n}^{[\mu n]}
\left( |\Delta \hat{f}(k)-\Delta \hat{f
}(-k)|\log k+\left| \Delta \hat{f}(k)\right| \log k\right) \right)  \nonumber \\
&&+O\left( \max\limits_{n\leq k\leq \lbrack \mu n]}|\hat{f}(k)|\log
k\right) .  \label{13}
\end{eqnarray}%
It is not difficult to see that
$$
\left\| \sum\limits_{k=n}^{[\mu n]-1}\hat{f}(-k-1)E_{k}(-x)\right\|
_{L}=I+O\left( n\max\limits_{n<k\leq \lbrack \mu
n]}|\hat{f}(-k)|\right) ,
$$
where
$$
I:=\int_{n^{-1}\leq |x|\leq \pi }\left| \frac{1}{2\sin (x/2)}%
\sum\limits_{k=n}^{[\mu
n]-1}\hat{f}(-k-1)e^{\frac{i(2k+1)x}{2}}\right| dx.
$$
Since the trigonometric function system is orthonormal, we have
\begin{eqnarray*}
I &\leq &\left( \int_{n^{-1}\leq |x|\leq \pi }\left|
\sum\limits_{k=n}^{[\mu
n]-1}\hat{f}(-k-1)e^{\frac{i(2k+1)x}{2}}\right| ^{2}dx\right) ^{1/2} \\
&&\hspace{3cm}\times \left( \int_{n^{-1}}^{\pi }\frac{1}{\sin ^{2}(x/2)}%
dx\right) ^{1/2} \\
&=&O\left( \sqrt{n}\left( \sum\limits_{k=n+1}^{[\mu n]}|\hat{f}%
(-k)|^{2}\right) ^{1/2}\right) \\
&=&O\left( n\max\limits_{n\leq k\leq \lbrack \mu
n]}|\hat{f}(-k)|\right) ,
\end{eqnarray*}%
which yields that
\begin{eqnarray}
\frac{1}{[\mu n]-n}\left\| \sum\limits_{k=n}^{[\mu n]-1}\hat{f}%
(-k-1)E_{k}(-x)\right\| _{L}=O\left( \max\limits_{n<k\leq \lbrack \mu n]}|%
\hat{f}(-k)|\right) .  \label{14}
\end{eqnarray}%
By combining (\ref{10}), (\ref{12}) - (\ref{14}), with Lemma 1 and
the condition
$$
\lim\limits_{n\rightarrow \infty }\Vert f-S_{n}(f)\Vert _{L}=0,
$$
and the fact (since $f\in L_{2\pi }$) that
$$
\lim\limits_{n\rightarrow \infty }\hat{f}(-n)=0,
$$
we have for $n\rightarrow \infty ,$
\begin{eqnarray}
\Vert \hat{f}(-n)E_{n}(-x)\Vert _{L} &\leq &\sum\limits_{k=n}^{[\mu
n]}|\Delta \hat{f}(k)-\Delta \hat{f}(-k)|\log k  \nonumber \\
&&+\Vert \tau _{\mu n,n}(f)-S_{n}(f)\Vert _{L}  \nonumber \\
&&+O\left( \max\limits_{[\lambda ^{-1}n]\leq k\leq \lbrack \lambda n]}|\hat{f%
}(k)|\log k\right)  \nonumber \\
&&+O\left( \max\limits_{n<k\leq \lbrack \mu n]}|\hat{f}(-k)|\right)
\nonumber
\\
&=&\sum\limits_{k=n}^{[\mu n]}|\Delta \hat{f}(k)-\Delta
\hat{f}(-k)|\log k+o(1).  \label{15}
\end{eqnarray}%
On the other hand, we have
\begin{eqnarray}
\Vert \hat{f}(-n)E_{n}(-x)\Vert _{L}\geq |\hat{f}(-n)|\Vert
D_{n}(x)\Vert _{L}\geq \frac{1}{\pi }|\hat{f}(-n)|\log n.
\label{16}
\end{eqnarray}%
Hence, from (\ref{15}), (\ref{16}), and (\ref{9}), we have that
$$
|\hat{f}(-n)|\log n\leq \sum\limits_{k=n}^{[\mu n]}|\Delta
\hat{f}(k)-\Delta \hat{f}(-k)|\log k\leq \varepsilon
$$
holds for sufficiently large $n$, which, together with (\ref{12}),
completes the proof of necessity.

In view of Lemma 1, we can see that the condition (\ref{2}) in
Theorem 1 can be replaced by the following condition
$$
\lim\limits_{\mu \rightarrow 1^{+}}\limsup\limits_{n\rightarrow
\infty }\sum\limits_{k=n}^{[\mu n]}|\Delta \hat{f}(-k)|\log k=0,
$$
and the proof of the result is easier. Therefore we have a corollary
to Theorem 1.

{\bf Corollary 1.}\quad {\it
Let $f(x)\in L_{2\pi }$ be a complex valued function. If both $\{\hat{f}%
(n)\}_{n=0}^{+\infty }\in MVBVS$ and $\{\hat{f}(-n)\}_{n=0}^{+\infty
}\in MVBVS$, then
$$
\lim\limits_{n\rightarrow \infty }\Vert f-S_{n}(f)\Vert _{L}=0
$$
if and only if
$$
\lim\limits_{n\rightarrow \infty }\hat{f}(n)\log |n|=0.
$$
}

If $f(x)$ is a real valued function, then its Fourier coefficients $\hat{f}%
(n)$ and $\hat{f}(-n)$ are a pair of conjugate complex numbers.
Consequently, $\{\hat{f}(n)\}_{n=0}^{+\infty }\in MVBVS$ if and only if $\{%
\hat{f}(-n)\}_{n=0}^{+\infty }\in MVBVS$. Thus, we have the
following generalization of the classical result (cf. Result Two in
the introduction):

{\bf Corollary 2.}\quad {\it Let $f(x)\in L_{2\pi }$ be a real
valued even function and (\ref{1a}) be its Fourier series. If
${\mathbf{A}}=\{a_{n}\}_{n=0}^{+\infty }\in MVBVS$ in real
sense, i.e. $\{a_{n}\}$ is a nonnegative sequence, and there is a number $%
\lambda \geq 2$ such that
$$
\sum\limits_{k=m}^{2m}|\Delta a_{k}|\leq C(\mathbf{A})\frac{1}{m}%
\sum\limits_{k=[\lambda ^{-1}m]}^{[\lambda m]}a_{k}
$$
for all $n=1,2,\ldots ,$ then
$$
\lim\limits_{n\rightarrow \infty }\Vert f-S_{n}(f)\Vert _{L}=0
$$
if and only if
$$
\lim\limits_{n\rightarrow \infty }a_{n}\log n=0.
$$
}

\section{$L^{1}$ Approximation}

Let $E_{n}(f)_{L}$ be the best approximation of a complex valued
function $f\in L_{2\pi }$ by trigonometric polynomials of degree $n$
in $L^{1}$ norm, that is,
$$
E_{n}(f)_{L}:=\inf\limits_{c_{k}}\left\|
f-\sum\limits_{k=-n}^{n}c_{k}e^{ikx}\right\| _{L}.
$$

We establish the corresponding $L^{1}-$approximation theorem in a
similar way to Theorem 1:

{\bf Theorem 2.}\quad {\it Let $f(x)\in L_{2\pi }$ be a complex
valued function, $\{\psi _{n}\}$ a decreasing sequence tending to
zero with
\begin{eqnarray}
\psi _{n}\sim \psi _{2n},  \label{17}
\end{eqnarray}%
i.e., there exist positive constants $C_{1}$ and $C_{2},$ such that $%
C_{1}\psi _{n}\leq \psi _{2n}\leq C_{2}\psi _{n}.$ If both $\{\hat{f}%
(n)\}_{n=0}^{+\infty }\in MVBVS$ and $\{\hat{f}(-n)\}_{n=0}^{+\infty
}\in MVBVS$, then
\begin{eqnarray}
\Vert f-S_{n}(f)\Vert _{L}=O(\psi _{n})  \label{18}
\end{eqnarray}%
if and only if
\begin{eqnarray}
E_{n}(f)_{L}=O(\psi _{n})\;\;\;\mbox{and}\;\;\;\hat{f}(n)\log
|n|=O(\psi _{|n|}).  \label{19}
\end{eqnarray}
}

{\bf Proof.} Under the condition of Theorem 2, we see from
(\ref{11a}) in the proof of Theorem 1 that
\begin{eqnarray*}
\Vert f-S_{n}(f)\Vert _{L} &\leq &\Vert f-\tau _{\mu n,n}(f)\Vert
_{L}+O\left( \max\limits_{[\lambda ^{-1}n]\leq |k|\leq \lbrack \lambda n]}|%
\hat{f}(k)|\log |k|\right) \\
&\leq &C(\mu )E_{n}(f)+O\left( \max\limits_{[\lambda ^{-1}n]\leq
|k|\leq \lbrack \lambda n]}|\hat{f}(k)|\log |k|\right) ,
\end{eqnarray*}%
thus (\ref{18}) holds if (\ref{17}) and (\ref{19}) hold. Now if
(\ref{18}) holds, then
$$
E_{n}(f)_{L}=O(\psi _{n})
$$
and
$$
\Vert f-\tau _{\mu n,n}(f)\Vert _{L}=O(\psi _{n}).
$$
From (\ref{6}) - (\ref{8}) in the proof of Lemma 4 and condition
(\ref{17}), we have
\begin{eqnarray}
|\hat{f}(n)|\log n &\leq &\frac{C(\lambda
)}{n}\sum\limits_{j=1}^{[\lambda n]-[\lambda ^{-1}n]+1}\left\|
f-S_{[\lambda ^{-1}n]+j}(f)\right\| _{L}
\label{20} \\
&&\;\;\;\;+C(\lambda )\Vert f-S_{n}(f)\Vert _{L}  \nonumber \\
&=&O(\psi _{n}).  \nonumber
\end{eqnarray}%
Since $\{\hat{f}(-n)\}_{n=0}^{+\infty }\in MVBVS,$ by a similar argument to (%
\ref{20}), we also have
$$
\left| \hat{f}(-n)\right| \log n=O(\psi _{n}).
$$
This completes the proof of Theorem 2.

In particular, if we take
$$
\psi _{n}:=\frac{1}{(n+1)^{r}}\omega \left(
f^{(r)},\frac{1}{n+1}\right) _{L},
$$
where $r$ is a positive integer, and $\omega (f,t)_{L}$ is the
modulus of continuity of $f$ in $L^{1}$ norm, i.e.
$$
\omega (f,t)_{L}:=\max\limits_{0\leq h\leq t}\Vert f(x+h)-f(x)\Vert
_{L}.
$$
By Theorem 2 and \ the Jackson theorem (e.g. see \cite{Xie-Zhou2} or \cite%
{Zygmund}) in $L^{1}-$space, we immediately have

{\bf Corollary 3.}\quad {\it
Let $f(x)\in L_{2\pi }$ be a complex valued function. If both $\{\hat{f}%
(n)\}_{n=0}^{+\infty }\in MVBVS$ and $\{\hat{f}(-n)\}_{n=0}^{+\infty
}\in MVBVS$ hold, then
$$
\Vert f-S_{n}(f)\Vert _{L}=O\left( \frac{1}{(n+1)^{r}}\omega \left( f^{(r)},%
\frac{1}{n+1}\right) _{L}\right)
$$
if and only if
$$
\hat{f}(n)\log |n|=O\left( \frac{1}{(n+1)^{r}}\omega \left( f^{(r)},\frac{1}{%
n+1}\right) _{L}\right) .
$$
}

This corollary generalizes the corresponding results in \cite{Xie-Zhou} and %
\cite{Le-Zhou}.

\

\begin{flushleft}
\rm Dan Sheng Yu and Song Ping Zhou:\\
Institute of Mathematics\\
Zhejiang  Sci-Tech University\\
Xiasha Economic Development Area\\
Hangzhou, Zhejiang 310018  China\\
And\\
Department of Mathematics, Statitics \& Computer Science\\
St. Francis Xavier University\\
Antigonish, Nova Scotia, Canada B2G 2W5\\
Email: dsyu@zjip.com (D. S. Yu)\\
\hspace{14mm}szhou@zjip.com (S. P. Zhou)
\end{flushleft}

\begin{flushleft}
\rm Ping Zhou:\\
 Department of Mathematics, Statitics \& Computer Science\\
St. Francis Xavier University\\
Antigonish, Nova Scotia, Canada B2G 2W5\\
Email: pzhou@stfx.ca
\end{flushleft}

\end{document}